\documentclass[reqno,11pt]{amsart}
\usepackage{amsmath, latexsym, amsfonts, amssymb, amsthm, amscd}

\numberwithin{equation}{section}

\newtheorem{theorem}{Theorem}
\newtheorem{lemma}[theorem]{Lemma}
\newtheorem{proposition}[theorem]{Proposition}

\newtheorem{rem}{Remark}

\newtheorem{assumption}{Assumption}

\newcommand{\bbL}{{\ensuremath{\mathbb L}} }

\newcommand{\bbP}{{\ensuremath{\mathbb P}} }

\newcommand{\cC}{{\ensuremath{\mathcal C}} }

\newcommand{\cG}{{\ensuremath{\mathcal G}} }

\newcommand{\ga}{\alpha}
\newcommand{\gb}{\beta}

\newcommand{\gd}{\delta}

\newcommand{\gep}{\varepsilon}       

\newcommand{\gr}{\rho}

\newcommand{\gs}{\sigma}

\newcommand{\gp}{\varphi}

\DeclareMathSymbol{\leqslant}{\mathalpha}{AMSa}{"36} 
\DeclareMathSymbol{\geqslant}{\mathalpha}{AMSa}{"3E} 
\DeclareMathSymbol{\eset}{\mathalpha}{AMSb}{"3F}     
\newcommand{\dd}{\text{\rm d}}             

\DeclareMathOperator*{\union}{\bigcup}       

\newcommand{\R}{\mathbb{R}}
\newcommand{\Z}{\mathbb{Z}}
\newcommand{\N}{\mathbb{N}}

\def\bP{\ensuremath{\bs{\mathrm{P}}}}
\def\bE{\ensuremath{\bs{\mathrm{E}}}}

\newcommand{\ind}{\bs{1}}

\def\bs{\boldsymbol}

\def\lf{\lfloor}
\def\rf{\rfloor}

\def\oT{\overline{T}}
\def\oH{\overline{H}}
\def\hG{\widehat{G}}
\def\tC{\mathtt{C}}

\newcommand{\mysim}{\overset{\star}{\sim}}

\title[LLT for random walks conditioned to stay positive]{A local limit theorem for random walks\\conditioned to stay positive}

\author{Francesco Caravenna}

\address{Universit\`a di Milano-Bicocca, Dipartimento di Matematica e Applicazioni, \mbox{Edificio} U5,\break via Cozzi 53, 20125 Milano, Italy \hfill \break
\indent \textit{and} \hfill\break
\indent Laboratoire de Probabilit{\'e}s de P 6 \& 7 
and  Universit{\'e} Paris 7
, U.F.R. Mathematiques, Case 7012, 2 place Jussieu, 75251 Paris cedex 05, France\hfill\break
\indent \textit{Home page:} {\tt http://www.matapp.unimib.it/\~{}fcaraven/}}

\email{f.caravenna\@@sns.it}



\begin{document}

\begin{abstract}
We consider a real random walk $S_n=X_1+\ldots +X_n$ attracted (without centering) to the normal law: this means that for a suitable norming sequence $a_n$ we have the weak convergence $S_n/a_n \Rightarrow \gp(x)\dd x$, $\gp(x)$ being the standard normal density. A local refinement of this convergence is provided by Gnedenko's and Stone's Local Limit Theorems, in the lattice and nonlattice case respectively.

Now let $\cC_n$ denote the event $(S_1 > 0, \ldots, S_n > 0 )$ and let $S_n^+$ denote the random variable $S_n$ conditioned on $\cC_n$: it is known that $S_n^+/a_n \Rightarrow \gp^+(x) \dd x$, where $\gp^+(x):=x \exp(-x^2/2)\ind_{(x\geq 0)}$. What we establish in this paper is an equivalent of Gnedenko's and Stone's Local Limit Theorems for this weak convergence. We also consider the particular case when $X_1$ has an absolutely continuous law: in this case the uniform convergence of the density of $S_n^+/a_n$ towards $\gp^+(x)$ holds under a standard additional hypothesis, in analogy to the classical case. We finally discuss an application of our main results to the asymptotic behavior of the joint renewal measure of the ladder variables process. Unlike the classical proofs of the LLT, we make no use of characteristic functions: our techniques are rather taken from the so--called Fluctuation Theory for random walks.
%
%
\end{abstract}

\maketitle

\section{Introduction and results} \label{sec:intro}

\subsection{The nonlattice case} Let $S_n=X_1+\ldots +X_n$ be a real random walk attracted (without centering) to the normal law. This means that $\{X_k\}$ is an IID sequence of real random variables, and for a suitable norming sequence $a_n$ we have the weak convergence
\begin{equation} \label{eq:clt}
    S_n/a_n \Rightarrow \gp(x)\,\dd x\,, \qquad \gp(x) := \frac{1}{\sqrt{2\pi}}\, e^{-x^2/2}\,.
\end{equation}
This is the case for example when $\bE(X_1)=0$ and $\bE(X_1^2) =: \gs^2 \in (0,\infty)$ with $a_n:= \gs \sqrt{n}$, by the Central Limit Theorem.

\smallskip

We recall that, by the standard theory of stability \cite[\S{}IX.8 \& \S{}XVII.5]{cf:Fel2}, for equation \eqref{eq:clt} to hold it is necessary and sufficient that $\bE(X_1)=0$, that the truncated variance $V(t):=\bE(X_1^2 \, \ind_{(|X_1|\leq t)})$ be \textsl{slowly varying} at~$\infty$ (that is $V(ct)/V(t) \to 1$ as $t\to\infty$ for every $c>0$) and that the sequence $a_n$ satisfy the condition $a_n^2 \sim n V(a_n)$ as $n\to\infty$.

\smallskip

For the moment we assume that the law of~$X_1$ is \textsl{nonlattice}, that is not supported in $(b+c\Z)$ for any $b\in\R, c>0$. Then a local refinement of \eqref{eq:clt} is provided by Stone's Local Limit Theorem \cite{cf:Sto65,cf:Sto67}, that in our notations reads as (cf.~\cite[\S{}8.4]{cf:BinGolTeu})
\begin{equation} \label{eq:llt}
    a_n\, \bP \big( S_n \in [x,x+h)\, \big) = h\, \gp(x/a_n) + o(1) \qquad (n\to\infty)\,,
\end{equation}
\textsl{uniformly} for $x\in\R$ and $h$ in compact sets in~$\R^+$.

\smallskip

In this paper we consider the asymptotic behavior of the random walk $\{S_n\}$ conditioned to stay positive. More precisely, let $\cC_n:=(S_1 > 0, \ldots, S_n > 0)$ and let $S_n^+$ denote the random variable $S_n$ under the conditional probability $\bP(\,\cdot\, |\,\cC_n)$: if \eqref{eq:clt} holds then one has an analogous weak convergence result for $S_n^+/a_n$, namely
\begin{equation} \label{eq:pos_clt}
    S_n^+/a_n \Rightarrow \gp^+(x)\,\dd x\,, \qquad \gp^+(x) := x\, e^{-x^2/2} \,\ind_{(x\geq 0)}\,.
\end{equation}
This is an immediate consequence of the fact \cite{cf:Igl74,cf:Bol76,cf:Don85} that, whenever \eqref{eq:clt} holds, the whole process $\{S_{\lf nt \rf} / a_n\}_{t\in[0,1]}$ under $\bP(\,\cdot\, |\,\cC_n)$ converges weakly as $n\to\infty$ to the standard Brownian meander process $\{B_t^+\}_{t\in[0,1]}$, and $\gp^+(x)\,\dd x$ is the law of $B_1^+$, cf.~\cite{cf:RevYor}.

\medskip

Our main result is an analogue of Stone's LLT for the weak convergence~\eqref{eq:pos_clt}.

\smallskip

\begin{theorem} \label{th:main}
If $X_1$ is nonlattice and \eqref{eq:clt} holds, then
\begin{equation}\label{eq:pos_llt}
a_n\, \bP \big( S_n \in [x,x+h)\, \big| \,\cC_n \big) = h\, \gp^+(x/a_n) + o(1) \qquad (n\to\infty)\,,
\end{equation}
uniformly for $x\in\R$ and $h$ in compact sets in~$\R^+$.
\end{theorem}

\smallskip

The main difficulty with respect to the classical case is given by the fact that under the conditional probability $\bP(\,\cdot\, |\,\cC_n)$ the increments of the walk $\{X_k\}$ are no longer independent. This is a major point in that the standard proof of Stone's LLT relies heavily on characteristic functions methods. As a matter of fact, we make no use of characteristic functions: our methods are rather of combinatorial nature, and we make an essential use of the so--called Fluctuation Theory for random walks. The core of our proof consists in expressing the law of $S_n$ under $\bP(\,\cdot\, |\,\cC_n)$ as a suitable mixture of the laws of $\{S_k\}_{0\leq k \leq n}$ \textsl{under the unconditioned measure}~$\bP$, to which Stone's LLT can be applied. Thus our ``Positive LLT'' is in a sense directly derived from Stone's LLT.

\smallskip

%

We point out that our methods may in principle be applied to the case when the random walk is attracted to a generic stable law (the analogue of~\eqref{eq:pos_clt} in this case is also provided by \cite{cf:Don85}), so that it should be possible to obtain an equivalent of Theorem~\ref{th:main} in this general setting.


\smallskip

\subsection{The lattice case}

Let us consider now the lattice case: we assume that $X_1$ is supported in $(b+c\Z)$, for the least such~$c$. In this case the local version of~\eqref{eq:clt} is given by Gnedenko's Local Limit Theorem \cite[\S{}8.4]{cf:BinGolTeu}, which says that
\begin{equation} \label{eq:lat_llt}
    \frac{a_n}{c}\, \bP\big( S_n=bn+cx \big) = \gp\big((bn+cx)/a_n\big) + o(1) \qquad (n\to\infty)\,,
\end{equation}
uniformly for $x\in\Z$.

\medskip

We can derive the local version of~\eqref{eq:pos_clt} also in this setting.

\begin{theorem} \label{th:lat_main}
If $X_1$ is lattice with span~1 and \eqref{eq:clt} holds, then
\begin{equation*}
\frac{a_n}{c}\, \bP\big( S_n=bn+cx \,\big|\, \cC_n \big) = \gp^+\big((bn+cx)/a_n\big) + o(1) \qquad (n\to\infty)\,,
\end{equation*}
uniformly for $x\in\Z$.
\end{theorem}

The proof is omitted since it can be recovered from the proof of Theorem~1 with only slight modifications (some steps are even simpler).


\smallskip

\subsection{The density case}

When the law of $X_1$ is absolutely continuous with respect to Lebesgue measure and \eqref{eq:clt} holds, one may ask whether the density of $S_n/a_n$ converges to $\gp(x)$ in some pointwise sense. However, it is easy to build examples \cite[\S{}46]{cf:GneKol} satisfying \eqref{eq:clt}, such that for every~$n$ the density of $S_n/a_n$ is unbounded in any neighborhood of~0: therefore without some extra--assumption one cannot hope for convergence to hold at each point. Nevertheless, if one looks for the \textsl{uniform convergence} of the density, then there is a simple condition which turns out to be necessary and sufficient.

\begin{assumption} \label{ass:main}
The law of $X_1$ is absolutely continuous, and for some $k\in\N$ the density $f_k(x)$ of $S_k$ is essentially bounded: $f_k(x) \in L^\infty(\R,\dd x)$.
\end{assumption}

It is easy to see that if this assumption holds, then for large~$n$ the density $f_n(x)$ admits a bounded and continuous version. A proof that Assumption~\ref{ass:main} yields the uniform convergence of the (continuous versions of the) density of $S_n/a_n$ towards~$\gp(x)$, namely
\begin{equation*}
    \sup_{x\in\R} \big| a_n f_n(a_n x) - \gp(x) \big| \to 0 \qquad (n\to\infty)\,,
\end{equation*}
can be found in~\cite[\S{}46]{cf:GneKol}. On the other side, the necessity of Assumption~\ref{ass:main} for the above convergence to hold is evident.

\medskip

We can derive a completely analogous result for $S_n^+$.

\begin{theorem} \label{th:den_main}
Assume that $X_1$ satisfies Assumption~\ref{ass:main}, and that \eqref{eq:clt} holds. Then:
\begin{enumerate}
\item $S_n^+$ has an absolutely continuous law, whose density $f_n^+(x)$ is bounded and continuous (except at $x=0$) for large~$n$;
\item the (continuous version of the) density of $S_n^+/a_n$ converges uniformly to~$\gp^+(x)$:
\begin{equation*}
    \sup_{x\in\R} \big| a_n f_n^+(a_n x) - \gp^+(x) \big| \to 0 \qquad (n\to\infty)\,.
\end{equation*}
\end{enumerate}
\end{theorem}

This Theorem can be proved following very closely the proof of Theorem~1: in fact equation \eqref{eq:main} in Section~\ref{sec:first} provides an explicit expression for $f_n^+(x)$, that can be shown to converge to $\gp^+(x)$ with the very same arguments given in Section~\ref{sec:second}.


\smallskip

\subsection{Asymptotic behavior of the ladder renewal measure}

As a by--product of the Local Limit Theorems described above, we have a result on the asymptotic behavior of the renewal measure of the ladder variables process. For simplicity we take the arithmetic setting, assuming that $X_1$ is supported by~$\Z$ and it is aperiodic, but everything works similarly in the general lattice and nonlattice cases. The renewal mass function $u(n,x)$ of the ladder variables process is defined for $n\in\N,\ x\in\Z$ by
\begin{equation} \label{eq:def_u}
    u(n,x) := \sum_{r=0}^\infty \bP \big( T_r=n, H_r=x \big) = \bP \big( n \text{ is a ladder epoch}, S_n = x \big)\,,
\end{equation}
where $\{(T_k,H_k)\}$ is the (strict, ascending) ladder variables process associated to the random walk (the definitions are given in Section~\ref{sec:fluct}). Generalizing some earlier result of \cite{cf:Kee}, in \cite{cf:AliDon} it has been shown that, for $\{x_n\}$ such that $x_n/a_n\to 0$,
\begin{equation} \label{eq:as_ren}
    u(n,x_n) \sim \frac{1}{\sqrt{2\pi}\, n\, a_n}\, U(x_n-1) \sim \frac 1n\, \bP\big( S_n=x_n \big) \, U(x_n-1) \qquad (n\to\infty)\,,
\end{equation}
where $U(x) := \sum_{r=0}^\infty \bP (H_r \leq x)$ is the distribution function of the renewal measure associated to the ladder heights process (as a matter of fact, the proof of \eqref{eq:as_ren} given in~\cite{cf:AliDon} is carried out under the assumption that $\bE(X_1^2)<\infty$, but it can be easily extended to the general case).

\medskip

With the methods of the present paper we are able to show that the same relation is valid for $x = O(a_n)$, with no further restriction on~$X_1$ other than the validity of~\eqref{eq:clt}.

\begin{theorem} \label{th:as_ren}
Let $X_1$ be arithmetic with span~1 and such that equation \eqref{eq:clt} holds. Then for $x\in\Z$
\begin{equation} \label{eq:ext_as_ren}
    u(n,x) = \frac 1n\, \bP\big( S_n=x \big) \, U(x-1) \, \big( 1+o(1) \big) \qquad (n\to\infty)\,,
\end{equation}
uniformly for $x/a_n \in [\gep, 1/\gep]$, for every fixed $\gep>0$.
\end{theorem}

The proof of this theorem is a direct consequence of Theorem~\ref{th:lat_main}: the details are worked out in Section~\ref{sec:as_ren}.

\smallskip

Notice that in the r.h.s. of \eqref{eq:ext_as_ren} we could as well write $U(x)$ instead of $U(x-1)$, since $x\to\infty$ as $n\to\infty$. Also observe that putting together equation~\eqref{eq:as_ren} with Theorem~\ref{th:as_ren} one has the stronger result that equation~\eqref{eq:ext_as_ren} holds uniformly for $x/a_n \in [0,K]$, for every fixed $K>0$.

\smallskip

We point out that Theorem~\ref{th:as_ren} has been obtained also in \cite{cf:BryDon}, where the authors study random walks conditioned to stay positive in a different sense.


\smallskip

\subsection{Outline of the paper}

The exposition is organized as follows: in Section~\ref{sec:fluct} we recall some basic facts on Fluctuation Theory and stable laws, and we set the relative notation; we also give the proof of Theorem~\ref{th:as_ren}. The rest of the paper is devoted to the proof of Theorem~\ref{th:main}, which has been split in two parts. The first one, in Section~\ref{sec:first}, contains the core of the proof: using Fluctuation Theory we obtain an alternative expression for the law of~$S_n^+$, see equation~\eqref{eq:main}, and we prove a crucial weak convergence result connected to the renewal measure of the ladder variables process. Then in Section~\ref{sec:second} we apply these preliminary results, together with Stone's LLT, to complete the proof. Finally, some minor points have been deferred to the appendix.


\bigskip

\section{Fluctuation Theory and some applications} \label{sec:fluct}

In this section we are going to recall some basic facts about Fluctuation Theory for random walks, especially in connection with the theory of stable laws, and to derive some preliminary results. Standard references on the subject are~\cite{cf:Fel2} and~\cite{cf:BinGolTeu}.


\smallskip

\subsection{Regular variation}
A positive sequence $d_n$ is said to be \textsl{regularly varying} of index $\ga\in\R$ (we denote this by $d_n \in R_{\ga}$) if $d_n \sim L_n\, n^\ga$ as $n\to \infty$, where $L_n$ is \textsl{slowly varying} at~$\infty$ in that $L_{\lf tn\rf}/L_n\to 1$ as $n\to\infty$, for every $t>0$. If $d_n \in R_\ga$ with $\ga\neq 0$, up to asymptotic equivalence we can (and will) always assume \cite[Th.1.5.3]{cf:BinGolTeu} that $d_n = d(n)$, with $d(\cdot)$ a continuous, strictly monotone function, whose inverse will be denoted by~$d^{-1}(\cdot)$. Observe that if $d_n \in R_\ga$ then $d^{-1}(n) \in R_{1/\ga}$ and $1/d_n \in R_{-\ga}$.

\smallskip

Let us recall two basic facts on regularly varying sequences that will be used a number of times in the sequel. The first one is a uniform convergence property \cite[Th.1.2.1]{cf:BinGolTeu}: if $d_n \in R_\ga$, then
\begin{equation}\label{eq:RV_unif_conv}
    d_{\lf tn \rf} = t^\ga \,d_n \,\big(1+o(1)\big) \qquad (n\to\infty)\,,
\end{equation}
uniformly for $t\in[\gep, 1/\gep]$, for every fixed $\gep>0$. The second basic fact \cite[Prop.1.5.8]{cf:BinGolTeu} is that if $d_n\in R_\ga$ with $\ga>-1$, then
\begin{equation}\label{eq:RV_asymp}
    \sum_{k=1}^n d_k \sim \frac{n d_n}{\ga+1} \qquad (n\to\infty)\,.
\end{equation}


\smallskip

\subsection{Ladder variables and stability}
The first (strict ascending) \textsl{ladder epoch} $T_1$ of a random walk $S_n = X_1 + \ldots + X_n$ is the first time the random walk enters the positive half line, and the corresponding \textsl{ladder height} $H_1$ is the position of the walk at that time:
\[
    T_1:=\inf \{n> 0: S_n >0\} \qquad H_1:=S_{T_1}\,.
\]
Iterating these definitions one gets the following ladder variables: more precisely, for $k>1$ one defines inductively
\[
    T_{k}:=\inf \{n>T_{k-1}: S_n >H_{k-1}\} \qquad H_k:=S_{T_k}\,,
\]
and for convenience we put $(T_0,H_0):=(0,0)$. The \textsl{weak} ascending ladder variables are defined in a similar way, just replacing $>$ by $\geq$ in the relations $(S_n > 0)$ and $(S_n > H_{k-1})$ above. In the following we will rather consider the weak \textsl{descending} ladder variables $(\oT_k, \oH_k)$, which are by definition the weak ascending ladder variables of the walk $\{-S_n\}$. Observe that, by the strong Markov property, both $\{(T_k, H_k)\}_k$ and $\{(\oT_k, \oH_k)\}_k$ are bidimensional renewal processes, that is random walks on $\R^2$ with step law supported in the first quadrant.

\smallskip

It is known that $X_1$ is in the domain of attraction (without centering) of a stable law if and only if $(T_1,H_1)$ lies in a bivariate domain of attraction, cf. \cite{cf:GreOmeTeu,cf:DonGre,cf:Don95}. This fact will play a fundamental role in our derivation: let us specialize it to our setting. By hypothesis $X_1$ is attracted to the normal law, that is $S_n/a_n \Rightarrow \gp(x) \,\dd x$, so that by the standard theory of stability $a_n \in R_{1/2}$. We define two sequences $b_n,\ c_n$ by
\begin{equation} \label{eq:def_b_c}
    \log \frac{n}{\sqrt 2} = \sum_{m=1}^\infty \frac{\gr_m}{m}\, e^{-\frac{m}{b_n}} \qquad c_n := a(b_n)\,,
\end{equation}
where $\gr_m := \bP(S_m > 0)$: then $b_n \in R_2$, $c_n \in R_1$ and we have the weak convergence
\begin{equation} \label{eq:weak_conv}
    \bigg( \frac{T_n}{b_n}\,, \frac{H_n}{c_n} \bigg) \Rightarrow Z\,, \qquad \bP\big(Z \in (\dd x,\,\dd y)\big) = \frac{e^{-1/2x}}{\sqrt{2\pi}\, x^{3/2}} \ind_{(x\geq 0)} \,\dd x\cdot \, \gd_1(\dd y)\,,
\end{equation}
where $\gd_1(\dd y)$ denotes the Dirac measure at $y=1$.

\smallskip

Thus the first ladder epoch $T_1$ is attracted to the positive stable law of index~$1/2$, as for the simple random walk case:
\begin{equation*}
    \frac{T_n}{b_n} \Rightarrow Y, \qquad \bP \big( Y \in \dd x \big) =  \frac{e^{-1/2x}}{\sqrt{2\pi}\, x^{3/2}} \ind_{(x\geq 0)} \,\dd x\,,
\end{equation*}
while for $\{H_k\}$ one has a generalized law of large numbers, with norming sequence~$c_n$: $H_n/c_n \Rightarrow 1$ (that is $H_1$ is \textsl{relatively stable}, cf.~\cite[\S{}8.8]{cf:BinGolTeu}).

\smallskip

We stress that we choose the sequence $a_n$ to be increasing, and by \eqref{eq:def_b_c} $b_n$ and $c_n$ are increasing too. We also recall that the norming sequence $b_n$ is sharply linked to the probability tail of the random variable~$T_1$, by the relation
\begin{equation}\label{eq:rel_stable}
    \bP\big(T_1 > b_n\big)\sim\sqrt{\frac{2}{\pi}}\; \frac{1}{n}\,.
\end{equation}
In fact, this condition is necessary and sufficient in order that a sequence $b_n$ be such that $T_n/b_n \Rightarrow Y$, cf.~\cite[\S{}XIII.6]{cf:Fel2}.

\begin{rem} \rm
It has already been noticed that when the step~$X_1$ has finite (nonzero) variance and zero mean,
\begin{equation*}
     \bE\big(X_1\big)=0 \qquad \bE\big(X_1^2\big)=:\gs^2 \in\, (0,\infty)\,,
\end{equation*}
by the Central Limit Theorem one can take $a_n = \gs\sqrt{n}$ in order that equation \eqref{eq:clt} holds. In other words, $X_1$ is in the \textsl{normal domain of attraction} of the normal law. In this case the first ladder height~$H_1$ is integrable~\cite{cf:Don80} and the behavior of the tail of~$T_1$ is given by
\begin{equation*}
    \bP\big(T_1 > n\big) \sim \frac{2\, \bE(H_1)}{\gs \sqrt{2\pi}} \frac{1}{\sqrt{n}} \qquad (n\to\infty)\,,
\end{equation*}
cf.~\cite[Th.1 in \S{}XII.7 \& Th.1 in \S{}XVIII.5]{cf:Fel2}.
This means that also $T_1$ and $H_1$ belong to the normal domain of attraction of their respective limit law, and one can take
\begin{equation*}
    b_n = \frac{\bE(H_1)^2}{\gs^2}\,n^2 \qquad c_n = \bE(H_1)\, n
\end{equation*}
in order that \eqref{eq:weak_conv} holds (we have used the law of large numbers for $H_1$ and relation \eqref{eq:rel_stable} for~$T_1$).
\end{rem}


\smallskip

\subsection{An asymptotic result}

As an application of the results exposed so far, we derive the asymptotic behavior of $\bP(\cC_n)$ as $n\to\infty$, which will be needed in the sequel. The connection with Fluctuation Theory is given by the fact that
\begin{equation*}
    \cC_n := \big(S_1>0, \ldots, S_n>0\big) = \big(\oT_1 > n\big)\,.
\end{equation*}
In analogy to what we have seen for~$T_1$, the fact that the random walk is attracted to the normal law implies that $\oT_1$ lies in the domain of attraction of the positive stable law of index~$1/2$. Therefore $\bP(\cC_n) \in R_{-1/2}$, and denoting by $\psi(t):=\bE(\exp (-t\oT_1))$ the Laplace transform of~$\oT_1$, by standard Tauberian theorems \cite[Ex.(c) in \S{}XIII.5]{cf:Fel2} we have that
\begin{equation*}
    \bP(\cC_n) \sim \frac{1}{\sqrt{\pi}}\, \big(1-\psi(1/n)\big) \qquad (n\to\infty)\,.
\end{equation*}
Now, for $\psi(t)$ we have the following explicit expression \cite[Th.1 in \S{}XII.7]{cf:Fel2}:
\begin{equation*}
    -\log (1-\psi(t)) \;=\; \sum_{m=1}^\infty \frac{\overline{\gr}_m}{m} e^{-mt} \;=\; -\log(1-e^{-t}) - \sum_{m=1}^\infty \frac{\gr_m}{m} e^{-mt}\,,
\end{equation*}
where $\overline{\gr}_m := \bP(S_m \leq 0)$. A look to \eqref{eq:def_b_c} then yields the desired asymptotic behavior:
\begin{equation}\label{eq:as_cn}
    \bP(\cC_n) \sim \frac{1}{\sqrt{2\pi}}\, \frac{b^{-1}(n)}{n} \qquad (n\to\infty)\,.
\end{equation}


\smallskip

\subsection{Two combinatorial identities}
The power of Fluctuation Theory for the study of random walks is linked to some fundamental identities. The most famous one is the so-called Duality Lemma \cite[\S{}XII]{cf:Fel2} which can be expressed as
\begin{equation} \label{eq:comb_id1}
\bP\big(n\ \text{is a ladder epoch}, S_n \in \dd x\big) = \bP\big(\cC_n, S_n \in \dd x\big) \,,
\end{equation}
where by $(n\ \text{is a ladder epoch})$ we mean of course the disjoint union $\cup_{k\geq 0}(T_k=n)$, and by $\bP(A,Z\in \dd x)$ we denote the finite measure $B\mapsto \bP(A,Z\in B)$. A second important identity, recently discovered by Alili and Doney \cite{cf:AliDon}, will play a fundamental role for us:
\begin{equation} \label{eq:comb_id2}
\bP \big(T_k=n, H_k \in \dd x\big) = \frac{k}{n}\, \bP \big(H_{k-1} < S_n \leq H_k, S_n \in \dd x\big) \,.
\end{equation}

We point out that both the above identities are of purely combinatorial nature, in the sense that they can be proved by relating the events on the two sides with suitable one to one, measure preserving transformations on the sample paths space.


\smallskip

\subsection{Proof of Theorem~\ref{th:as_ren}} \label{sec:as_ren}

We recall that by hypothesis $\gep$ is a fixed positive number. We start from the definition \eqref{eq:def_u} of $u(n,x)$: applying the Duality Lemma~\eqref{eq:comb_id1} we get
\begin{equation} \label{eq:u_1}
    u(n,x) = \bP \big(\,\cC_n, S_n = x\big) = \bP \big(\,\cC_n\big) \,  \bP \big( S_n = x \,\big|\, \cC_n \big)\,.
\end{equation}
Observe that
\begin{equation*}
    \inf_{z\in[\gep,1/\gep]} \gp^+(z) > 0 \qquad \inf_{z\in[\gep,1/\gep]} \gp(z) > 0\,,
\end{equation*}
which implies that both Theorem~\ref{th:lat_main} and Gnedenko's LLT~\eqref{eq:lat_llt} hold also in a ratio sense, namely
\begin{align*}
    \bP\big( S_n=x \,\big|\, \cC_n \big) &= \frac{1}{a_n}\, \gp^+(x/a_n) \, \big( 1+o(1) \big) \qquad (n\to\infty)\\
    \bP\big( S_n=x  \,\big) &= \frac{1}{a_n}\, \gp\,(x/a_n) \, \big( 1+o(1) \big) \qquad (n\to\infty)\,,
\end{align*}
uniformly for $x/a_n\in[\gep, 1/\gep]$. Since $\gp^+(z) =\sqrt{2\pi}\, z\, \gp(z)$ for $z>0$, it follows that
\begin{equation} \label{eq:u_2}
    \bP\big( S_n=x \,\big|\, \cC_n \big) = \sqrt{2\pi}\, \frac{x}{a_n}\, \bP\big( S_n=x  \,\big) \, \big( 1+o(1) \big)  \qquad (n\to\infty)\,,
\end{equation}
uniformly for $x/a_n\in[\gep, 1/\gep]$.

The asymptotic behavior of $\bP(\cC_n)$ is given by~\eqref{eq:as_cn}, and comparing equation \eqref{eq:ext_as_ren} with \eqref{eq:u_2} and \eqref{eq:u_1} we are left with proving that
\begin{equation*}
    U(x) = x\, \frac{b^{-1}(n)}{a(n)} \, \big( 1+o(1) \big) \qquad (n\to\infty)\,,
\end{equation*}
uniformly for $x/a_n\in[\gep, 1/\gep]$. We recall that $U(x)$ is the distribution function of the renewal measure associated to the ladder height process $\{H_k\}$, which is relatively stable, since $H_n/c_n \Rightarrow 1$ as $n\to\infty$. Then Theorem~8.8.1 in \cite{cf:BinGolTeu} gives that $U(x) \sim c^{-1}(x)$ as $x\to\infty$, hence it rests to show that
\begin{equation*}
    \frac{x}{c^{-1}(x)} \, \frac{b^{-1}(n)}{a(n)} \to 1 \qquad (n\to\infty)\,,
\end{equation*}
uniformly for $x/a_n\in[\gep, 1/\gep]$, or equivalently, setting $x=z\,a_n$, that
\begin{equation*}
    \frac{z\,b^{-1}(n)}{c^{-1}(z\,a(n))} \to 1 \qquad (n\to\infty)\,,
\end{equation*}
uniformly for $z\in[\gep,1/\gep]$. However, as $c^{-1}(\cdot)\in R_1$, by \eqref{eq:RV_unif_conv} we have that
\begin{equation*}
    c^{-1}(z\,a(n)) \sim z\,c^{-1}(a(n)) \qquad (n\to\infty)\,,
\end{equation*}
uniformly for $z\in[\gep,1/\gep]$, and the proof is completed observing that $c^{-1}(a(n)) = b^{-1}(n)$, by the definition \eqref{eq:def_b_c} of~$c_n$.\qed


\bigskip

\section{Proof of Theorem~\ref{th:main}: first part} \label{sec:first}

\subsection{A fundamental expression}

We are going to use Fluctuation Theory to express the law of $S_n^+$ in a more useful way. For $x>0$ and $n>1$ we have
\begin{align*}
    &n \,\bP \big( \cC_n,\, S_n \in \dd x \big) \overset{\eqref{eq:comb_id1}}{=} n \,\bP \big( n \text{ is a ladder epoch},\, S_n \in \dd x \big)\\
     &\qquad = \sum_{r=1}^\infty n\,\bP \big( T_r=n,\, S_n \in \dd x \big) \overset{\eqref{eq:comb_id2}}{=} \sum_{r=1}^\infty r\,\bP \big( H_{r-1} < x \leq H_r,\, S_n \in \dd x \big)\,,
\end{align*}
where we have used both the combinatorial identities \eqref{eq:comb_id1}, \eqref{eq:comb_id2}. With a simple manipulation we get
\begin{align*}
    & \sum_{r=1}^\infty r\,\bP \big( H_{r-1} < x \leq H_r,\, S_n \in \dd x \big) = \sum_{r=1}^\infty \sum_{k=0}^{r-1} \bP \big( H_{r-1} < x \leq H_r,\, S_n \in \dd x \big) \\
    &\qquad = \sum_{k=0}^\infty \sum_{r=k+1}^{\infty} \bP \big( H_{r-1} < x \leq H_r,\, S_n \in \dd x \big) = \sum_{k=0}^\infty \bP \big( H_{k} < x,\, S_n \in \dd x \big)\,,
\end{align*}
and using the Markov property
\begin{align*}
    \bP \big( H_{k} < x,\, S_n \in \dd x \big) = \sum_{m=0}^{n-1} \int_{[0,x)} \bP \big( T_k=m,\, H_{k} \in \dd z\big) \, \bP \big(S_{n-m} \in \dd x - z \big)\,.
\end{align*}
In conclusion we obtain the following relation (which is essentially the same as equation~(10) in~\cite{cf:AliDon}):
\begin{align}
    & \bP \big( S_n/a_n \in \dd x\,\big|\, \cC_n \big) \nonumber\\
    &\qquad =\; \frac{1}{n \bP(\cC_n)} \, \sum_{m=0}^{n-1} \int_{[0,a_n x)} \Bigg( \sum_{k=0}^\infty \bP \big( T_k=m,\, H_{k} \in \dd z\big) \Bigg) \, \bP \big(S_{n-m} \in a_n\dd x - z \big) \nonumber\\
    &\qquad =\; \frac{b^{-1}(n)}{n \bP(\cC_n)} \: \int_{[0,1)\times[0,x)} \dd\mu_n(\ga,\gb)\; \bP \bigg( \frac{S_{\lf n(1-\ga) \rf}}{a_n} \in \dd x - \gb \bigg) \label{eq:main} \,,
\end{align}
where $\mu_n$ is the finite measure on $[0,1)\times [0,\infty)$ defined by
\begin{equation}\label{eq:def_mun}
    \mu_n (A) \;:=\; \frac{1}{b^{-1}(n)} \, \sum_{k=0}^\infty \bP \bigg( \bigg( \frac{T_k}{n}, \frac{H_k}{a_n} \bigg) \in A \bigg)\,,
\end{equation}
for $n\in\N$ and for any Borel set $A\subseteq [0,1) \times [0,\infty)$. Notice that $\mu_n$ is nothing but a suitable rescaling of the renewal measure associated to the ladder variables process. Also observe that the sum defining~$\mu_n$ can be stopped at $k=n-1$, since by definition $T_k \geq k$ for every~$k$; hence $\mu_n$ is indeed a finite measure.

\smallskip

Before proceeding, we would like to stress the importance of equation~\eqref{eq:main}, which is in a sense the core of our proof. The reason is that in the r.h.s. the conditioning on~$\cC_n$ has disappeared: we are left with a mixture, governed by the measure $\mu_n$, of the laws of $\{S_{\lf n(1- \ga) \rf}\}_{\ga\in[0,1)}$ \textsl{without conditioning}, and the asymptotic behavior of these laws can be controlled with Stone's Local Limit Theorem \eqref{eq:llt} (if we exclude the values of~$\ga$ close to~1).

\smallskip

In the following subsection we study the asymptotic behavior of the sequence of measures $\{\mu_n\}$, and in the next section we put together these preliminary results to conclude the proof of Theorem~\ref{th:main}.


\smallskip

\subsection{A weak convergence result} We are going to show that as $n\to\infty$ the sequence of measure $\{\mu_n\}$ converges weakly to the finite measure $\mu$ defined by
\begin{equation}\label{eq:def_mu}
    \mu(A) \;:=\; \int_A \dd\ga\,\dd\gb\, \frac{\gb}{\sqrt{2\pi}\,\ga^{3/2}} \, e^{-\gb^2/2\ga}\,,
\end{equation}
for any Borel set $A\subseteq [0,1)\times[0,\infty)$ (it is easy to check that $\mu$ is really a finite measure, see below). Since we are not dealing with probability measures, we must be most precise: we mean weak convergence with respect to the class $C_b$ of bounded and continuous functions on~$\R^2$: $\mu_n \Rightarrow \mu$ iff $\int h \,\dd\mu_n \to \int h \,\dd\mu$ for every $h\in C_b$. If we introduce the distribution functions $F_n,\, F$ of the measures $\mu_n,\, \mu$:
\begin{equation*}
    F_n(a,b) := \mu_n \big( [0,a] \times [0,b] \big) \qquad F(a,b) := \mu \big( [0,a] \times [0,b] \big)\,,
\end{equation*}
then proving that $\mu_n \Rightarrow \mu$ as $n\to\infty$ is equivalent to showing that $F_n(a,b) \to F(a,b)$ for every $(a,b)\in[0,1]\times[0,\infty]$ (notice that $\infty$ is included, because the total mass of~$\mu_n$ is not fixed).

\begin{proposition} \label{prop:weak_conv}
The sequence of measures $\{\mu_n\}$ converges weakly to the measure~$\mu$.
\end{proposition}

\proof
We start checking the convergence of the total mass:
\begin{equation*}
    F_n(1,\infty) = \frac{1}{b^{-1}(n)} \, \sum_{k=0}^\infty \bP \big( T_k \leq n \big) =: \frac{1}{b^{-1}(n)} \, G(n)\,,
\end{equation*}
where $G(n)$ is the distribution function of the renewal measure associated to the ladder epochs process~$\{T_k\}$. There is a sharp link between the asymptotic behavior as $n\to\infty$ of $G(n)$ and that of $\bP(T_1>n)$, given by \cite[Lem. in \S{}XIV.3]{cf:Fel2}:
\begin{equation} \label{eq:ren_meas}
    G(n) \sim \frac{2}{\pi} \, \frac{1}{\bP (T_1 > n)} \qquad (n\to\infty)\,.
\end{equation}
Since from relation \eqref{eq:rel_stable} we have that
\begin{equation*}
    \bP \big(T_1 > n\big) \sim \sqrt{\frac{2}{\pi}} \, \frac{1}{b^{-1}(n)} \qquad (n\to\infty)\,,
\end{equation*}
it follows that $F_n(1,\infty) \to \sqrt{2/\pi}$ as $n\to\infty$. On the other hand, the check that $F(1,\infty) = \sqrt{2/\pi}$ is immediate:
\begin{equation*}
    F(1,\infty) = \frac{1}{\sqrt{2\pi}} \int_0^1 \dd\ga\, \frac{1}{\ga^{3/2}} \int_0^\infty \dd\gb \,\gb\,e^{-\gb^2/2\ga}  = \frac{1}{\sqrt{2\pi}} \int_0^1 \dd\ga\, \frac{1}{\sqrt{\ga}} = \sqrt{\frac{2}{\pi}}\,.
\end{equation*}

\smallskip

Since the total mass converges, we claim that it suffices to show that
\begin{equation} \label{eq:claim}
    \liminf_{n\to\infty} \mu_n \big( (a_1,a_2] \times (b_1,b_2] \big) \;\geq\; \mu \big( (a_1,a_2] \times (b_1,b_2] \big)
\end{equation}
for all $0<a_1<a_2<1$, $0<b_1<b_2<\infty$, and weak convergence will be proved. The (simple) proof of this claim can be found in Appendix~\ref{app:claim}.

\smallskip

Directly from the definition of~$\mu_n$ we have
\begin{align*}
    \mu_n\big( (a_1,a_2] \times (b_1,b_2] \big) = \frac{1}{b^{-1}(n)} \sum_{k=0}^\infty \bP \bigg( \frac{T_k}{n} \in (a_1, a_2],\, \frac{H_k}{a_n} \in (b_1,b_2]  \bigg)\,.
\end{align*}
We simply restrict the sum to the set of $k$ such that $k/b^{-1}(n)\in (b_1+\gep, b_2-\gep]$, $\gep$ being a small fixed positive number, getting
\begin{equation} \label{eq:riemann}
    \mu_n\big( (a_1,a_2] \times (b_1,b_2] \big) \;\geq\; \frac{1}{b^{-1}(n)} \sum_{s\in\frac{\Z}{b^{-1}(n)} \cap (b_1+\gep, b_2-\gep]} \xi_n(s)\,,
\end{equation}
where
\begin{equation*}
    \xi_n(s) := \bP \bigg( \frac{T_{\lf s b^{-1}(n)\rf}}{n} \in (a_1, a_2],\, \frac{H_{\lf s b^{-1}(n)\rf}}{a_n} \in (b_1,b_2]  \bigg) \,.
\end{equation*}
By the definition \eqref{eq:def_b_c} of $c_n$, we have that $a_n = c(b^{-1}(n))$: then, using the weak convergence \eqref{eq:weak_conv} and the uniform convergence property of regularly varying sequences \eqref{eq:RV_unif_conv}, it is not difficult to check that
\begin{equation*}
    \xi_n(s) \to \bP \bigg( Y \in \bigg( \frac{a_1}{s^2}, \frac{a_2}{s^2}\bigg] \bigg) =: \xi(s) \qquad (n\to\infty)\,,
\end{equation*}
\textsl{uniformly} for $s\in (b_1+\gep, b_2-\gep]$.

Observe that the term in the r.h.s. of \eqref{eq:riemann} is a Riemann sum of the function~$\xi_n(s)$ over the bounded interval $(b_1+\gep, b_2-\gep]$. Since the sequence of functions $\{\xi_n(s)\}$ is clearly equibounded and converges uniformly to~$\xi(s)$, it is immediate to check that the r.h.s. of~\eqref{eq:riemann} does converge to the integral of~$\xi(s)$ over $(b_1+\gep, b_2-\gep]$. Therefore
\begin{align*}
    &\liminf_{n\to\infty} \mu_n\big( (a_1,a_2] \times (b_1,b_2] \big) \;\geq\; \int_{b_1+\gep}^{b_2-\gep} \dd s \; \bP \bigg( Y \in \bigg( \frac{a_1}{s^2}, \frac{a_2}{s^2}\bigg] \bigg)\\
    &\qquad =\; \int_{b_1+\gep}^{b_2-\gep} \dd s \int_{a_1/s^2}^{a_2/s^2} \dd z\; \frac{e^{-1/2z}}{\sqrt{2\pi}\,z^{3/2}} \;=\; \int_{b_1+\gep}^{b_2-\gep} \dd s \int_{a_1}^{a_2} \dd t\; \frac{s\,e^{-s^2/2t}}{\sqrt{2\pi}\,t^{3/2}}\\
    &\qquad =\; \mu \big( (a_1,a_2] \times (b_1+\gep, b_2-\gep] \big)\,,
\end{align*}
and letting $\gep\to 0$ relation \eqref{eq:claim} follows.\qed


\bigskip

\section{Proof of Theorem~\ref{th:main}: second part} \label{sec:second}


\subsection{General strategy}

Now we are ready to put together the results obtained in the last section. We start by rephrasing relation~\eqref{eq:pos_llt}, which is our final goal, in terms of $S_n/a_n$, a form that is more convenient for our purposes: we have to prove that
\begin{equation}\label{eq:goal}
    \forall K>0\quad \limsup_{n\to\infty}\, a_n \bigg[\; \sup_{x\in\R^+,\, h \leq K/a_n} \Big|\, \bP \big( S_n/a_n \in x+I_h\, \big| \,\cC_n \big) \;-\; h\, \gp^+(x) \,\Big| \; \bigg] = 0\,,
\end{equation}
where $I_h:=[0,h)$, and $x+I_h := [x,x+h)$.

\medskip

Altough the idea behind the proof is quite simple, our arguments depend on an approximation parameter~$\gep$ and there are a number of somewhat technical points. In order to keep the exposition as transparent as possible, it is convenient to introduce the following notation: given two real functions $f(n,x,h,\gep)$ and $g(n,x,h,\gep)$ of the variables $n\in\N$, $x\in\R^+$, $h\in\R^+$ and $\gep\in(0,1)$, we say that $f \mysim g$ if and only if
\begin{equation*}
     \forall K>0\quad \limsup_{\gep\to 0}\, \limsup_{n\to\infty}\, a_n \bigg[ \sup_{x\in\R^+,\,h\leq K/a_n} \big| f(n,x,h,\gep) - g(n,x,h,\gep) \big| \bigg] = 0\,.
\end{equation*}
With this terminology we can reformulate \eqref{eq:goal} as
\begin{equation} \label{eq:goal1}
    \bP \big( S_n/a_n \in x+I_h\, \big| \,\cC_n \big) \;\mysim\; h\:\gp^+(x)\,.
\end{equation}

\medskip

To obtain a more explicit expression of the l.h.s. of \eqref{eq:goal1}, we resort to equation~\eqref{eq:main}: with an easy integration we get
\begin{equation} \label{eq:start}
\bP \big( S_n/a_n \in x+I_h\, \big| \,\cC_n \big) \;=\; \frac{b^{-1}(n)}{n \bP(\cC_n)} \: \int_{D_1^{x+h}} \dd\mu_n(\ga,\gb)\; \hG^{x,h}_n(\ga,\gb) \,,
\end{equation}
where we have introduced the notation $D_a^b := [0,a)\times[0,b)$, and
\begin{equation} \label{eq:def_Ghat}
    \hG^{x,h}_n(\ga,\gb) \;:=\; \bP \bigg( \frac{S_{\lf n(1-\ga)\rf}}{a_n} \in \big\{ (x - \gb)+ I_h \big\} \cap [0,\infty) \bigg) \,.
\end{equation}
In order to determine the asymptotic behavior of the r.h.s. of \eqref{eq:start}, we recall that:
\begin{itemize}
\item from~\eqref{eq:as_cn} we have
\begin{equation*}
    \frac{b^{-1}(n)}{n \bP(\cC_n)} \to \sqrt{2\pi}\,;
\end{equation*}
\item from Proposition~\ref{prop:weak_conv} we have that $\mu_n \Rightarrow \mu$;
\item from Stone's LLT \eqref{eq:llt} it follows that, for large~$n$, $\hG^{x,h}_n(\ga,\gb)$ is close to
\begin{equation} \label{eq:def_G}
    \cG^{x,h}(\ga,\gb) \;:=\; h \, \frac{1}{\sqrt{1-\ga}}\; \gp \bigg( \frac{x-\gb}{\sqrt{1-\ga}} \bigg)\,,
\end{equation}
where we have used that $a_{n(1-\ga)} \sim \sqrt{1-\ga}\; a_n$ as $n\to\infty$, by \eqref{eq:RV_unif_conv}.
\end{itemize}
In fact, the rest of this section is devoted to showing that
\begin{equation} \label{eq:diff_goal}
    \bP \big( S_n/a_n \in x+I_h\, \big| \,\cC_n \big) \;\mysim\; \sqrt{2\pi} \int_{D_1^{x}} \dd\mu(\ga,\gb)\; \cG^{x,h}(\ga,\gb) \,.
\end{equation}
It may not be a priori obvious whether this coincides with our goal \eqref{eq:goal1}, that is whether
\begin{equation}\label{eq:conj}
    \gp^+(x) \;=\; \sqrt{2\pi} \int_{D_1^x} \dd\mu(\ga,\gb)\; \frac{1}{\sqrt{1-\ga}}\; \gp \bigg( \frac{x-\gb}{\sqrt{1-\ga}} \bigg)\,.
\end{equation}
Indeed this relation holds true: in fact \eqref{eq:diff_goal} implies the weak convergence of $S_n/a_n$ under $\bP(\,\cdot\, |\, \cC_n)$ towards a limiting law with the r.h.s. of~\eqref{eq:conj} as density, and we already know from \eqref{eq:pos_clt} that $S_n/a_n$ under $\bP(\,\cdot\, |\, \cC_n)$ converges weakly to $\gp^+(x) \,\dd x$. Anyway, a more direct verification of \eqref{eq:conj} is also given in Appendix~\ref{app:integral}.

\medskip

Thus we are left with proving \eqref{eq:diff_goal}, or equivalently
\begin{equation*}
    \int_{D_1^{x+h}} \dd\mu_n(\ga,\gb)\; \hG^{x,h}_n(\ga,\gb) \;\mysim\; \int_{D_1^x} \dd\mu(\ga,\gb)\; \cG^{x,h}(\ga,\gb)\,.
\end{equation*}
Since $\mysim$ is an equivalence relation, this will be done through a sequence of intermediate equivalences:
\begin{equation*}
    \int_{D_1^{x+h}} \dd\mu_n\; \hG^{x,h}_n  \;\mysim\; \ldots \;\mysim\; \ldots \;\mysim\; \ldots \;\mysim\; \int_{D_1^{x}} \dd\mu\; \cG^{x,h}\,,
\end{equation*}
and for ease of exposition the proof has been accordingly split in four steps. The idea is quite simple: we first restrict the domain from $D_1^{x+h}$ to $D_{1-\gep}^x$ (steps~1--2), then we will be able to apply Stone's LLT and Proposition~\ref{prop:weak_conv} to pass from $(\hG_n^{x,h}, \mu_n)$ to $(\cG^{x,h},\mu)$ (step~3), and finally we come back to the domain~$D_1^x$ (step~4).

\smallskip

Before proceeding, we define a slight variant $G_n^{x,h}$ of $\hG_n^{x,h}$:
\begin{equation}\label{eq:def_Gn}
    G^{x,h}_n(\ga,\gb) \;:=\; \bP \bigg( \frac{S_{\lf n(1-\ga)\rf}}{a_n} \in (x - \gb)+ I_h \bigg)
\end{equation}
(notice that we have simply removed the set $[0,\infty)$, see~\eqref{eq:def_Ghat}) and we establish a preliminary lemma.

\begin{lemma} \label{lem:prel}
For every $K>0$ there exists a positive constant $\tC = \tC(K)$ such that
\begin{equation*}
    G^{x,h}_n (\ga,\gb) \leq \frac{\tC}{a_{\lf(1-\ga)n\rf}} \qquad \forall n\in\N,\ \forall x,\gb\in\R,\ \forall \ga\in[0,1),\ \forall h \leq K/a_n\,,
\end{equation*}
and the same relation holds also for~$\hG_n^{x,h}(\ga,\gb)$.
\end{lemma}

\proof
Since by definition $\hG_n^{x,h}(\ga,\gb) \leq G_n^{x,h}(\ga,\gb)$, it suffices to prove the relation for~$G_n^{x,h}$. However, this is a simple consequence of Stone's LLT \eqref{eq:llt}, that we can rewrite in terms of $S_n/a_n$ as
\begin{equation} \label{eq:llt_alt}
    \forall K>0\quad \limsup_{l\to\infty}\, a_l \bigg[\; \sup_{y\in\R,\, h' \leq K/a_l} \Big|\, \bP \big( S_l/a_l \in y+I_{h'}\, \big) \;-\; h'\, \gp(y) \,\Big| \; \bigg] = 0\,.
\end{equation}
In fact from this relation, using the triangle inequality and the fact that $\sup_{x\in\R} |\gp(x)| < \infty$, it follows easily that for every~$K>0$
\begin{equation} \label{eq:lem_step}
    a_l\, \bP \big( S_l/a_l \in y+I_{h'}\, \big) \leq \tC \qquad \forall l\in\N,\ \forall y\in\R,\ \forall h' \leq K/a_l\,,
\end{equation}
for some positive constant $\tC=\tC(K)$. Now it suffices to observe that $G_n^{x,h}$ can be written as
\begin{equation} \label{eq:def_G_alt}
    G^{x,h}_n(\ga,\gb) \;=\; \bP \bigg( \frac{S_{\lf n(1-\ga)\rf}}{a_{\lf n(1-\ga)\rf}} \in \frac{a_n}{a_{\lf n(1-\ga)\rf}} (x - \gb)+ I_{\frac{h\,a_n}{a_{\lf n(1-\ga)\rf}}} \bigg)\,,
\end{equation}
so that we can apply \eqref{eq:lem_step} with $l=\lf n(1-\ga) \rf$ and analogous substitutions.\qed

\smallskip


\smallskip

\subsection{First step}

In the first intermediate equivalence we pass from the domain $D_1^{x+h}$ to $D_{1-\gep}^{x+h}$, that is we are going to show that
\begin{equation*}
    \int_{D_1^{x+h}} \dd\mu_n\; \hG^{x,h}_n  \;\mysim\; \int_{D_{1-\gep}^{x+h}} \dd\mu_n\; \hG^{x,h}_n \,.
\end{equation*}
This means by definition that for every $K>0$
\begin{equation}\label{eq:step1_alt}
    \limsup_{\gep\to 0} \limsup_{n\to\infty}\, R_n^\gep = 0 \,,
\end{equation}
where $R_n^\gep := \sup_{\{x\in\R^+, \; h\leq K/a_n\}} r_n^\gep(x,h)$ and
\begin{equation*}
    r_n^\gep(x,h) := a_n \int_{[1-\gep,1)\times[0,x+h)} \dd\mu_n(\ga,\gb)\; \hG^{x,h}_n (\ga,\gb)\,.
\end{equation*}

\smallskip

Applying Lemma~\ref{lem:prel} and enlarging the domain of integration, we get
\begin{equation} \label{eq:intermed}
\begin{split}
    R_n^\gep &\leq \tC\,a_n \int_{[1-\gep,1)\times[0,\infty)} \dd\mu_n(\ga,\gb)\; \frac{1}{a_{\lf(1-\ga)n\rf}} \\
    &= \tC\,a_n \sum_{m=\lf (1-\gep)n\rf}^{n-1} \Bigg[ \frac{1}{b^{-1}(n)} \, \sum_{k=0}^\infty \bP \big( T_k = m \big) \Bigg] \frac{1}{a_{n-m}} \\
    &= \tC\, \frac{a_n}{b^{-1}(n)} \sum_{m=\lf (1-\gep)n\rf}^{n-1} \frac{u(m)}{a_{n-m}}\,,
\end{split}
\end{equation}
where in the second line we have applied the definition \eqref{eq:def_mun} of~$\mu_n$, and in the third line we have introduced $u(m) := \sum_{k=0}^\infty \bP(T_k=m)$, which is the mass function of the renewal measure associated to the ladder epochs process~$\{T_k\}$. In the proof of Proposition~\ref{prop:weak_conv} we have encountered the asymptotic behavior of the distribution function $G(n):=\sum_{m=1}^n u(m)$, see~\eqref{eq:ren_meas}. The corresponding local asymptotic behavior for $u(m)$ follows since the sequence $u(m)$ is decreasing in~$m$ (this is a simple consequence of the Duality Lemma \eqref{eq:comb_id1}, see also~\cite[Th.4]{cf:Don97}): hence
\begin{equation*}
    u(m) \sim \frac{1}{\pi} \, \frac{1}{m\bP (T_1 > m)} \sim \frac{1}{\sqrt{2\pi}} \,\frac{b^{-1}(m)}{m} \qquad (m\to\infty)\,,
\end{equation*}
having used~\eqref{eq:rel_stable}. It follows that $u(m) \leq C_1\, b^{-1}(m) /m$ for every~$m$, for some positive constant~$C_1$. Recalling that $b^{-1}(\cdot)$ is increasing, from \eqref{eq:intermed} we get
\begin{align*}
    R_n^\gep & \;\leq\;  \tC\,C_1 \frac{a_n}{b^{-1}(n)} \sum_{m=\lf (1-\gep)n\rf}^{n-1} \frac{b^{-1}(m)}{m\,a_{n-m}}\\
    &\;\leq\; \tC C_1\, \frac{a_n}{\lf(1-\gep)n\rf} \sum_{k=1}^{\lf \gep n\rf} \frac{1}{a_k} \;\leq\; \tC C_1C_2 \,\frac{\gep}{1-\gep}\,\frac{a_n}{a_{\lf \gep n\rf}}  \,,
\end{align*}
for some positive constant $C_2$: in the last inequality we have used \eqref{eq:RV_asymp}, since $a_n \in R_{1/2}$. Now from \eqref{eq:RV_unif_conv} we have that $a_n/a_{\lf \gep n\rf} \to 1/\sqrt{\gep}$ as $n\to\infty$, hence
\begin{equation*}
    \limsup_{n\to\infty}\, R_n^\gep \leq C\, \frac{\sqrt{\gep}}{1-\gep}\,,
\end{equation*}
with~$C:=\tC C_1 C_2$, and \eqref{eq:step1_alt} follows.


\smallskip

\subsection{Second step}

Now we show that we can restrict the domain from $D_{1-\gep}^{x+h}$ to $D_{1-\gep}^{x}$:
\begin{equation*}
    \int_{D_{1-\gep}^{x+h}} \dd\mu_n\; \hG^{x,h}_n  \;\mysim\; \int_{D_{1-\gep}^x} \dd\mu_n\; \hG^{x,h}_n \;=\; \int_{D_{1-\gep}^x} \dd\mu_n\; G^{x,h}_n \,,
\end{equation*}
where the equality simply follows from the fact that by definition (see \eqref{eq:def_Ghat} and \eqref{eq:def_Gn})
\begin{equation*}
    \hG^{x,h}_n(\ga,\gb) = G^{x,h}_n(\ga,\gb) \qquad \text{for } \gb \leq x\,.
\end{equation*}
We have to show that for every $K>0$
\begin{equation}\label{eq:step2_alt}
    \limsup_{\gep\to 0} \,\limsup_{n\to\infty}\, Q_n^\gep = 0\,,
\end{equation}
where $Q_n^\gep := \sup_{\{x\in\R^+, \; h\leq K/a_n\}} q_n^\gep(x,h)$ and
\begin{equation*}
    q_n^\gep(x,h) := a_n \int_{[0,1-\gep)\times[x,x+h)} \dd\mu_n(\ga,\gb)\; \hG^{x,h}_n (\ga,\gb)\,.
\end{equation*}

\smallskip

From Lemma~\ref{lem:prel} and from the fact that $a_n$ is increasing we easily get
\begin{equation*}
    q_n^\gep(x,h) \leq \tC\; \frac{a_n}{a_{\lf\gep n\rf}} \; \mu_n \big( [0,1-\gep)\times[x,x+h) \big)\,.
\end{equation*}
As $a_n\in R_{1/2}$, we have $a_n/a_{\lf\gep n\rf}\to 1/\sqrt{\gep}$ as $n\to\infty$ by \eqref{eq:RV_unif_conv}, hence for fixed $\gep>0$ we can find a positive constant $C_1=C_1(\gep)$ such that for all~$n\in\N$
\begin{equation*}
    q_n^\gep(x,h) \leq \tC\, C_1 \; \mu_n \big( [0,1-\gep)\times[x,x+h) \big)\,.
\end{equation*}
However the term in the r.h.s. can be easily estimated: using the definition \eqref{eq:def_mun} of~$\mu_n$, for $h\leq K/a_n$ we get
\begin{align*}
    &\mu_n \big( [0,1-\gep)\times[x,x+h) \big) \;=\; \frac{1}{b^{-1}(n)} \sum_{k=0}^\infty \bP \big( T_k < (1-\gep)n,\, H_k \in [a_n x, a_n x + a_n h) \big) \\
    &\qquad \;\leq\; \frac{1}{b^{-1}(n)} \sum_{k=0}^\infty \bP \big( H_k \in [a_n x, a_n x + K) \big) \;\leq\; \frac{1}{b^{-1}(n)}\; \sup_{z\in\R^+} U\big( [z,z+K) \big)\,,
\end{align*}
where $U(\dd x) := \sum_{k=0}^\infty \bP (H_k \in \dd x)$ is the renewal measure associated to the ladder heights process $\{H_k\}$, that we have already encountered in the proof of Theorem~\ref{th:as_ren}. Notice that
\begin{equation*}
    \forall K>0 \qquad \sup_{z\in\R^+} U \big( [z,z+K)\big) =: C_2 < \infty \,,
\end{equation*}
which holds whenever $\{H_k\}$ is a transient random walk, cf.~\cite[Th.1 in \S{}VI.10]{cf:Fel2}. Thus for every fixed~$\gep > 0$
\begin{equation*}
    Q_n^\gep = \sup_{x\in\R^+,\, h\leq K/a_n} q_n^\gep(x,h) \;\leq\; \tC C_1 C_2\, \frac{1}{b^{-1}(n)} \;\to\; 0 \qquad (n\to\infty)\,,
\end{equation*}
and \eqref{eq:step2_alt} follows.


\smallskip

\subsection{Third step}

This is the central step: we prove that
\begin{equation*}
    \int_{D_{1-\gep}^x} \dd\mu_n\; G^{x,h}_n \;\mysim\; \int_{D_{1-\gep}^x} \dd\mu\; \cG^{x,h}\,,
\end{equation*}
that is for every $K>0$
\begin{equation}\label{eq:step3_alt}
    \limsup_{\gep\to 0} \, \limsup_{n\to\infty} \, \sup_{x\in\R^+,\, h\leq K/a_n} a_n\,\Bigg| \int_{D_{1-\gep}^x} \dd\mu_n\, G^{x,h}_n \;-\; \int_{D_{1-\gep}^x} \dd\mu\, \cG^{x,h} \Bigg| = 0\,.
\end{equation}

By the triangle inequality
\begin{equation}\label{eq:int1}
\begin{split}
    &a_n\,\Bigg| \int_{D_{1-\gep}^x} \dd\mu_n\, G^{x,h}_n \;-\; \int_{D_{1-\gep}^x} \dd\mu\, \cG^{x,h} \Bigg| \\
    &\qquad \leq\; a_n\,\int_{D_{1-\gep}^x} \dd\mu_n\, \big| G^{x,h}_n - \cG^{x,h} \big| \;+\; a_n\,\Bigg| \int_{D_{1-\gep}^x} \dd\mu_n\, \cG^{x,h} \;-\; \int_{D_{1-\gep}^x} \dd\mu\, \cG^{x,h} \Bigg|\,,
\end{split}
\end{equation}
and we study separately the two terms in the r.h.s. above.

\smallskip
\subsubsection{First term}

With a rough estimate we have
\begin{equation} \label{eq:int1.5}
\begin{split}
    &a_n\,\int_{D_{1-\gep}^x} \dd\mu_n\, \big| G^{x,h}_n - \cG^{x,h} \big|\\
    &\qquad \leq\; \bigg[ \sup_{n\in\N} \mu_n\big(D_{1}^\infty\big) \bigg] \, \Bigg( \sup_{(\ga,\gb)\in D_{1-\gep}^\infty} a_n\, \Big| G^{x,h}_n(\ga,\gb) - \cG^{x,h}(\ga,\gb) \Big| \Bigg)\,,
\end{split}
\end{equation}
and notice the prefactor in the r.h.s. is bounded since $\mu_n (D_1^\infty) \to \mu (D_1^\infty)$. For the remaining term, we use the triangle inequality and the definition \eqref{eq:def_G} of~$\cG^{x,h}$, getting
\begin{equation}\label{eq:int2}
\begin{split}
    &a_n\,\Big| G^{x,h}_n(\ga,\gb) - \cG^{x,h}(\ga,\gb) \Big|\\
    &\qquad\quad \;\leq\; \bigg(\frac{a_n}{a_{\lf(1-\ga)n\rf}}\bigg)\, a_{\lf(1-\ga)n\rf}\, \bigg| G^{x,h}_n(\ga,\gb) - \frac{h\,a_n}{a_{\lf(1-\ga)n\rf}} \, \gp \bigg(  \frac{a_n\, (x-\gb)}{a_{\lf(1-\ga)n\rf}} \bigg)  \bigg|\\
    &\qquad\qquad\ \ +\; (h\,a_n)\, \bigg| \frac{a_n}{a_{\lf(1-\ga)n\rf}} \, \gp \bigg(  \frac{a_n\, (x-\gb)}{a_{\lf(1-\ga)n\rf}} \bigg) - \frac{1}{\sqrt{1-\ga}} \, \gp \bigg( \frac{x-\gb}{\sqrt{1-\ga}} \bigg) \bigg|\,.
\end{split}
\end{equation}

Let us look at the first term in the r.h.s. above: by the by the uniform convergence property of regularly varying sequences \eqref{eq:RV_unif_conv} we have
\begin{equation} \label{eq:int3}
    \sup_{\ga\in(0,1-\gep)} \bigg| \frac{a_n}{a_{\lf(1-\ga)n\rf}} - \frac{1}{\sqrt{1-\ga}}\, \bigg| \to 0 \qquad (n\to\infty)\,,
\end{equation}
hence the prefactor is uniformly bounded. For the remaining part, from the expression \eqref{eq:def_G_alt} for~$G^{x,h}_n$ it is clear that one can apply Stone's LLT, see \eqref{eq:llt_alt}, yielding
\begin{equation*}
    \sup_{(\ga,\gb) \in D_{1-\gep}^\infty,\, x\in\R^+,\, h\leq K/a_n} \; a_{\lf(1-\ga)n\rf}\, \bigg| G^{x,h}_n(\ga,\gb) - \frac{h\,a_n}{a_{\lf(1-\ga)n\rf}} \, \gp \bigg(  \frac{a_n\, (x-\gb)}{a_{\lf(1-\ga)n\rf}} \bigg)  \bigg| \;\to\;0
\end{equation*}
as $n\to\infty$.

For the second term in the r.h.s. of \eqref{eq:int2}, notice that the prefactor $(h\, a_n)$ gives no problem since $h\leq K/a_n$ in our limit. On the other hand, it is easily seen that the absolute value is vanishing as $n\to\infty$, uniformly for $(\ga,\gb) \in D_{1-\gep}^\infty$ and for $x\in\R^+$: this is thanks to relation~\eqref{eq:int3} and to the fact that the function $\gp(x)$ is uniformly continuous. Coming back to equation \eqref{eq:int1.5}, we have shown that
\begin{equation}\label{eq:int4}
    \limsup_{n\to\infty} \, \sup_{x\in\R^+,\,h\leq K/a_n} a_n\,\int_{D_{1-\gep}^x} \dd\mu_n\, \big| G^{x,h}_n - \cG^{x,h} \big| \;=\; 0\,.
\end{equation}

\smallskip
\subsubsection{Second term} Using the definition \eqref{eq:def_G} of~$\cG^{x,h}$, the second term in the r.h.s. of equation \eqref{eq:int1} can be written as
\begin{equation}\label{eq:int5}
\begin{split}
    &a_n\,\Bigg| \int_{D_{1-\gep}^x} \dd\mu_n\, \cG^{x,h} \;-\; \int_{D_{1-\gep}^x} \dd\mu\, \cG^{x,h} \Bigg| \\
    &\qquad =\; (h\,a_n)\, \Bigg| \int_{D_{1-\gep}^\infty} \dd\mu_n\, \Psi(\ga, x-\gb) \;-\; \int_{D_{1-\gep}^\infty} \dd\mu\, \Psi(\ga, x-\gb)\; \Bigg|
\end{split}
\end{equation}
where we have introduced the shorthand
\begin{equation*}
    \Psi(s, t) := \frac{1}{\sqrt{1-s}} \, \gp \bigg( \frac{t}{\sqrt{1-s}} \bigg) \,\ind_{(t\geq 0)}\,.
\end{equation*}

As usual, for us $(h\,a_n)\leq K$ and we can thus concentrate on the absolute value in the r.h.s. of~\eqref{eq:int5}. Observe that, for fixed~$x\geq 0$, the function $(\ga,\gb) \mapsto \Psi(\ga, x-\gb)$ on the domain $D_{1-\gep}^\infty$ is bounded, and continuous except on the line $\gb=x$: since $\mu_n \Rightarrow \mu$, it follows that for fixed~$x$ the r.h.s. of~\eqref{eq:int5} is vanishing as $n\to\infty$. However, we would like the convergence to be uniform in~$x\in\R^+$: this stronger result holds true too, as one can verify by approximating $\Psi$ with a sequence of uniformly continuous functions (the details are carried out in Appendix~\ref{app:unif}). The net result is
\begin{equation} \label{eq:int6}
    \limsup_{n\to\infty} \, \sup_{x\in\R^+,\,h\leq K/a_n} \, a_n\,\Bigg| \int_{D_{1-\gep}^x} \dd\mu_n\, \cG^{x,h} \;-\; \int_{D_{1-\gep}^x} \dd\mu\, \cG^{x,h} \Bigg| \;=\; 0\,.
\end{equation}

\smallskip

Putting together relations \eqref{eq:int1}, \eqref{eq:int4} and \eqref{eq:int6} it is easily seen that \eqref{eq:step3_alt} holds (even without taking the limit in~$\gep$), and the step is completed.


\smallskip

\subsection{Fourth step}

We finally show that
\begin{equation*}
    \int_{D_{1-\gep}^{x}} \dd\mu\; \cG^{x,h} \;\mysim\; \int_{D_{1}^x} \dd\mu\; \cG^{x,h} \,,
\end{equation*}
that is, for every~$K>0$
\begin{equation} \label{eq:step4_alt}
    \limsup_{\gep\to 0}\, \limsup_{n\to\infty}\, \sup_{x\in\R^+,\, h\leq K/a_n} a_n \int_{[1-\gep,1)\times [0,x)}  \dd\mu (\ga,\gb)\, \cG^{x,h}(\ga,\gb) = 0\,.
\end{equation}

\smallskip

This is very easy: observe that
\begin{equation*}
    \cG^{x,h}(\ga,\gb) \leq \frac{h}{\sqrt{2\pi}\, \sqrt{1-\ga}}\,,
\end{equation*}
as one can check from the explicit expressions for $\cG^{x,h}$~\eqref{eq:def_G} and $\gp(x)$~\eqref{eq:clt}. Hence
\begin{equation*}
    a_n \int_{[1-\gep,1)\times [0,x)}  \dd\mu (\ga,\gb)\, \cG^{x,h}(\ga,\gb) \;\leq\; \frac{(h a_n)}{\sqrt{2\pi}} \int_{[1-\gep,1)\times [0,\infty)}  \dd\mu (\ga,\gb) \frac{1}{\sqrt{1-\ga}} \,,
\end{equation*}
and \eqref{eq:step4_alt} follows, because the function
\begin{equation*}
    \big\{ (\ga,\gb) \mapsto (1-\ga)^{-1/2} \big\} \in \bbL^1\big( D_1^\infty,\,\dd\mu \big)\,,
\end{equation*}
as on can easily verify. This completes the proof of Theorem~\ref{th:main}.


\appendix

\bigskip

\section{An elementary fact} \label{app:claim}

We prove the claim stated in the proof of Proposition~\ref{prop:weak_conv}, in a slightly more general context. Namely, let $\mu_n,\,\mu$ be finite measures on the domain $D:=[0,1)\times[0,\infty)$, with $\mu(\partial D)=0$. Assume that $\mu_n(D) \to \mu(D)$ as $n\to\infty$, and that
\begin{equation} \label{eq:claim1}
    \liminf_{n\to\infty} \mu_n \big( (a_1,a_2] \times (b_1,b_2] \big) \;\geq\; \mu \big( (a_1,a_2] \times (b_1,b_2] \big)\,,
\end{equation}
for all $0<a_1<a_2<1$, $0<b_1<b_2<\infty$. What we are going to show is that
\begin{equation} \label{eq:claim2}
    \exists\ \lim_{n\to\infty} \mu_n \big( (a_1,a_2] \times (b_1,b_2] \big) \;=\; \mu \big( (a_1,a_2] \times (b_1,b_2] \big)\,,
\end{equation}
for all $0<a_1<a_2<1$, $0<b_1<b_2<\infty$, and this implies that $\mu_n \Rightarrow \mu$.

\smallskip

Suppose that \eqref{eq:claim2} does not hold: then for some rectangle $Q:= (x_1,x_2] \times (y_1,y_2]$ contained in the interior of~$D$ and for some $\gep>0$ one has
\begin{equation} \label{eq:absurd}
    \limsup_{n\to\infty} \mu_n(Q) \;\geq\; \mu(Q) + \gep\,.
\end{equation}
We introduce for $\eta \in (0,1/2)$ the rectangle $W:=(\eta,1-\eta] \times (\eta, 1/\eta]$: by choosing~$\eta$ sufficiently small we can assume that $W \supseteq Q$ and that
\begin{equation}\label{eq:app_step1}
    \mu(W) \geq \mu(D)-\gep/2
\end{equation}
(we recall that by hypothesis $\mu(\partial D)=0$). The rectangle $W$ can be easily written as a \textsl{disjoint} union
\begin{equation*}
    W = Q \cup \union_{i=1}^4 Q_i\,,
\end{equation*}
where the rectangles $Q_i$ (whose exact definition however is immaterial) are defined by
\begin{eqnarray*}
    &Q_1 := (\eta,1-\eta]\times(\eta,y_1] \qquad &Q_2 := (\eta, x_1] \times (y_1,y_2]\\
    &Q_3 := (x_2, 1-\eta] \times (y_1,y_2] \qquad &Q_4 := (\eta,1-\eta]\times(y_2,1/\eta]\,.
\end{eqnarray*}
Now, on the one hand we have
\begin{equation*}
    \limsup_{n\to\infty} \mu_n(W) \;\leq\; \limsup_{n\to\infty} \mu_n(D) \;=\; \mu(D)\,,
\end{equation*}
but on the other hand
\begin{align*}
    &\limsup_{n\to\infty} \mu_n(W) \;=\; \limsup_{n\to\infty} \mu_n\bigg(Q \cup \union_{i=1}^4 Q_i \bigg) \;\geq\; \limsup_{n\to\infty} \mu_n(Q) + \liminf_{n\to\infty} \mu_n\bigg(\union_{i=1}^4 Q_i \bigg) \\
    &\qquad \overset{\eqref{eq:absurd}}{\geq} \mu(Q) + \gep + \sum_{i=1}^4 \liminf_{n\to\infty} \mu_n( Q_i) \;\overset{\eqref{eq:claim1}}{\geq} \mu(Q) + \gep + \sum_{i=1}^4 \mu(Q_i) \;=\; \gep + \mu(W) \\
    &\qquad \overset{\eqref{eq:app_step1}}{\geq} \mu(D) + \gep/2\,,
\end{align*}
which evidently is absurd, hence \eqref{eq:claim2} holds true.


\bigskip

\section{An integral} \label{app:integral}

We are going to give a more direct proof of relation~\eqref{eq:conj}: substituting the explicit expressions for $\gp(x)$, $\gp^+(x)$, $\mu$ given in equations \eqref{eq:clt}, \eqref{eq:pos_clt}, \eqref{eq:def_mu} and performing an elementary change of variable, we can rewrite it as
\begin{equation} \label{eq:conj1}
    x \, e^{-x^2/2} \;=\; \frac{x^2}{\sqrt{2\pi}} \int_0^1 \dd w \int_0^1 \dd z \; \frac{w}{z^{3/2}(1-z)^{1/2}} \, e^{-\frac{x^2}{2} \big[ \frac{w^2}{z} + \frac{(1-w)^2}{(1-z)} \big]}\,.
\end{equation}

Altough it is possible to perform explicitly the integration in the r.h.s. above, it is easier to proceed in a different way. Let $\{B_t\}$ be a standard Brownian motion and let $T_a := \inf\{t: B_t=a\}$ be its first passage time: then the law of $T_a$ is given by
\begin{equation*}
    \bbP \big( T_a \in \dd t \big) = g(a,t)\, \dd t\,, \qquad g(a,t) := \frac{a}{\sqrt{2\pi}\, t^{3/2}} e^{-a^2/2t}\,.
\end{equation*}
By the strong Markov property, for $x>0$ and $w\in(0,1)$ we have the equality in law $T_x \sim T_{wx} + T_{(1-w)x}$\,, where we mean that $T_{wx}$ and $T_{(1-w)x}$ are independent. Therefore
\begin{equation*}
    g(x,1) = \int_0^1 \dd z\, g\big(wx,z\big)\, g\big((1-w)x,1-z\big)\,,
\end{equation*}
and integrating over $w \in (0,1)$ we get
\begin{equation} \label{eq:app_rel1}
    g(x,1) = \int_0^1 \dd w \int_0^1 \dd z\, g\big(wx,z\big)\, g\big((1-w)x,1-z\big)\,.
\end{equation}

Now observe that relation \eqref{eq:conj1} can be written as
\begin{align*}
    g(x,1) &= \int_0^1 \dd w \int_0^1 \dd z\, \frac{1-z}{1-w}\, g\big(wx,z\big)\, g\big((1-w)x,1-z\big)\\
    &= \int_0^1 \dd w \int_0^1 \dd z\, \frac{z}{w}\, g\big(wx,z\big)\, g\big((1-w)x,1-z\big)\,,
\end{align*}
and comparing with \eqref{eq:app_rel1} we are left with showing that
\begin{equation*}
    \int_0^1 \dd w \int_0^1 \dd z\, \bigg( 1- \frac{z}{w} \bigg)\, g\big(wx,z\big)\, g\big((1-w)x,1-z\big)\;=\; 0 \\
\end{equation*}
However, the l.h.s. above can be decomposed in
\begin{equation*}
    \int_0^1 \dd w \int_w^1 \dd z\; \big( \ldots \big) \;+\; \int_0^1 \dd w \int_0^w \dd z\; \big( \ldots \big) \;=:\; I_1 \;+\; I_2 \,,
\end{equation*}
and with a change of variable one easily verifies that $I_1=-I_2$.



\bigskip

\section{A uniformity result} \label{app:unif}

We are going to show that
\begin{equation}\label{eq:app_but}
    \limsup_{n\to\infty} \, \sup_{x\in\R^+} \, \Bigg| \int_{D_{1-\gep}^\infty} \dd\mu_n\, \Psi(\ga, x-\gb) \;-\; \int_{D_{1-\gep}^\infty} \dd\mu\, \Psi(\ga, x-\gb)\; \Bigg| \;=\; 0\,,
\end{equation}
where we recall that $D_a^b := [0,a) \times [0,b)$ and the function $\Psi$ is defined by
\begin{equation*}
    \Psi(s, t) := \frac{1}{\sqrt{1-s}} \, \gp \bigg( \frac{t}{\sqrt{1-s}} \bigg) \,\ind_{(t\geq 0)}\,.
\end{equation*}

\smallskip

Let us consider the fixed domain $T:=[0,1-\gep]\times \R$. Here the function $\Psi$ is bounded, $\|\Psi\|_{\infty,T} = 1/\sqrt{2\pi\gep}$, and continuous except on the line $t=0$. We can easily build a family of approximations $\{\Psi_\gd\}$ of~$\Psi$ that are bounded and uniformly continuous on the whole~$T$, setting for $\gd>0$
\begin{equation*}
    \Psi_\gd (s,t) := \begin{cases}
\Psi(s,t) & t \geq 0\\
\Psi(s,0)\cdot (1+t/\gd) & t \in [-\gd, 0]\\
0 & t \leq -\gd
\end{cases} \,.
\end{equation*}
Notice that $\|\Psi_\gd\|_{\infty,T} = \|\Psi\|_{\infty,T}$, and that for $(s,t)\in T$
\begin{equation} \label{eq:app_int}
    \big|\Psi(s,t) - \Psi_\gd (s,t)\big| \leq \|\Psi\|_{\infty,T} \, \ind_{[-\gd, 0]}(t)\,.
\end{equation}

Let us introduce for short the notation $\Psi^x(\ga,\gb):=\Psi (\ga,x-\gb)$, and analogously for~$\Psi_\gd$. From the triangle inequality we get
\begin{equation} \label{eq:app_step2}
\begin{split}
    \Bigg| \int_{D_{1-\gep}^\infty} \dd\mu_n\, \Psi^x & \;-\; \int_{D_{1-\gep}^\infty} \dd\mu\, \Psi^x\; \Bigg| \;\leq\; \int_{D_{1-\gep}^\infty} \dd\mu_n \, \big| \Psi^x - \Psi^x_\gd \big| \\
    &\;+\; \int_{D_{1-\gep}^\infty} \dd\mu \, \big| \Psi^x - \Psi^x_\gd \big| \;+\; \Bigg| \int_{D_{1-\gep}^\infty} \dd\mu_n\, \Psi^x_\gd \;-\; \int_{D_{1-\gep}^\infty} \dd\mu\, \Psi^x_\gd\; \Bigg|\,.
\end{split}
\end{equation}
Using relation \eqref{eq:app_int}, the first two terms in the r.h.s. above can be estimated by
\begin{equation*}
    \|\Psi\|_{\infty,T} \, \Big( \mu_n \big( [0,1-\gep]\times[x,x+\gd] \big) + \mu \big( [0,1-\gep]\times[x,x+\gd] \big) \Big)\,.
\end{equation*}
Since $\mu$ is an absolutely continuous and finite measure, its distribution function is uniformly continuous: therefore for every $\eta>0$ we can take $\gd_0$ sufficiently small so that
\begin{equation*}
    \sup_{x\in\R^+}\, \mu \big( [0,1-\gep]\times[x,x+\gd_0] \big) \leq \frac{\eta}{4\|\Psi\|_{\infty,T}}\,.
\end{equation*}
On the other hand, we know that for every $x\geq 0$
\begin{equation*}
    \mu_n \big( [0,1-\gep]\times[x,x+\gd_0] \big) \to \mu \big( [0,1-\gep]\times[x,x+\gd_0] \big) \qquad (n\to\infty)\,,
\end{equation*}
and this convergence is uniform for $x\in\R^+$, as it can be easily checked. Hence by the triangle inequality we can choose $n_0$ so large that
\begin{equation*}
    \sup_{n\geq n_0} \, \sup_{x\in\R^+}\, \mu_n \big( [0,1-\gep]\times[x,x+\gd_0] \big) \leq \frac{\eta}{2\|\Psi\|_{\infty,T}}\,.
\end{equation*}
Finally, observe that for fixed $\gd_0$ the family of functions $\{\Psi^x_{\gd_0}\}_{x\in\R^+}$ is equibounded and equicontinuous: since $\mu_n \Rightarrow \mu$, from a classical result \cite[Cor. in \S{}VIII.1]{cf:Fel2} we have that the third term in the r.h.s. of \eqref{eq:app_step2} with $\gd=\gd_0$ is vanishing as $n\to\infty$ uniformly for $x\in\R^+$. Therefore we can assume that $n_0$ has been chosen so large that
\begin{equation*}
    \sup_{n\geq n_0} \,  \sup_{x\in\R^+}\, \Bigg| \int_{D_{1-\gep}^\infty} \dd\mu_n \, \Psi^x_{\gd_0} \;-\; \int_{D_{1-\gep}^\infty} \dd\mu\, \Psi^x_{\gd_0}\; \Bigg| \;\leq\; \frac{\eta}{4}\,.
\end{equation*}

Applying the preceding bounds to equation \eqref{eq:app_step2} with $\gd=\gd_0$, we have shown that for every $\eta>0$ we can find $n_0$ such that for every $n\geq n_0$
\begin{equation*}
    \sup_{x\in\R^+}\, \Bigg| \int_{D_{1-\gep}^\infty} \dd\mu_n\, \Psi^x \;-\; \int_{D_{1-\gep}^\infty} \dd\mu\, \Psi^x\; \Bigg| \;\leq\; \eta\,,
\end{equation*}
and equation \eqref{eq:app_but} is proved.


\bigskip

\section*{Acknowledgments}

I'm grateful to my Ph.D. supervisor Giambattista Giacomin for his help and suggestions and to Paolo Lorenzoni and Lorenzo Zambotti for their comments on the manuscript. After the submission of this paper I received reference~\cite{cf:BryDon} from Ron A. Doney and I would like to thank him for that. I also want to express all my gratitude to the referee for a number of important remarks that led to the improvement of several results.


\end{document}